\documentclass{aims}
\usepackage{amsmath}
  \usepackage{paralist}
 \usepackage[colorlinks=true]{hyperref}
\hypersetup{urlcolor=blue, citecolor=red}

\def\currentvolume{X}
 \def\currentissue{X}
  \def\currentyear{200X}
   \def\currentmonth{XX}
    \def\ppages{X--XX}

\usepackage{amssymb,amscd,amsthm,amsxtra}
\usepackage{dsfont,latexsym}
\usepackage{url}

\usepackage{mathrsfs}
\renewcommand{\mathcal}{\mathscr}

\definecolor{citation}{rgb}{0.2,0.5,0.2}
\definecolor{formula}{rgb}{0.1,0.2,0.5}
\definecolor{url}{rgb}{0,0.2,0.7}

\newtheorem{theorem}{Theorem}[section]
\newtheorem{corollary}[theorem]{Corollary}

\newtheorem{prop}[theorem]{Proposition}
\newtheorem{example}[theorem]{Example}

\theoremstyle{definition}

\theoremstyle{remark}

\numberwithin{equation}{section}

\newcommand{\CC}{\mathcal{C}}

\newcommand{\R}{{\mathds R}}
\newcommand{\N}{{\mathds N}}

\renewcommand{\le}{\leqslant}
\renewcommand{\leq}{\leqslant}
\renewcommand{\ge}{\geqslant}

\newcommand{\aand}{\,{\mbox{ and }}\,}
\newcommand{\dist}{\mbox{dist}}

\newlength{\defbaselineskip}
\newcommand{\setlinespacing}[1]
           {\setlength{\baselineskip}{#1 \defbaselineskip}}

\title
[Asymptotics of the $s$-perimeter as $s\searrow 0$]
{Asymptotics of the $s$-perimeter as $s\searrow 0$}

\author[S. Dipierro, A. Figalli, G. Palatucci and E. Valdinoci]{}

\subjclass{Primary: 49Q15, 35R09, 35R11;
Secondary: 45K05.}
\keywords{Nonlinear problems, nonlocal perimeter, fractional Laplacian, 
fractional Sobolev spaces, minimal surfaces}

\email{\href{mailto:dipierro@sissa.it}{dipierro@sissa.it}}
\email{\href{mailto:figalli@math.utexas.edu}{figalli@math.utexas.edu}}
\email{\href{mailto:giampiero.palatucci@unimes.fr}{giampiero.palatucci@unimes.fr}}
\email{\href{mailto:enrico@math.utexas.edu}{enrico@math.utexas.edu}}

\thanks{The first and the fourth author have been supported by FIRB 
``Project Analysis and Beyond''. The second author has been supported by
\href{http://www.nsf.gov/awardsearch/showAward.do?AwardNumber=0969962}{NSF
Grant~DMS-0969962}. The third author has been supported by
\href{http://prmat.math.unipr.it/~rivista/eventi/2010/ERC-VP/}{ERC grant
207573 ``Vectorial Problems''}. The fourth author has been supported by
ERC grant 277749 ``EPSILON Elliptic Pde's and Symmetry of Interfaces and
Layers for Odd Nonlinearities''.}

\begin{document}
\maketitle

\centerline{\scshape Serena Dipierro }
\medskip
{\footnotesize
 \centerline{SISSA - International School for Advanced Studies}
   \centerline{Sector of Mathematical Analysis}
   \centerline{Via Bonomea, 265}
   \centerline{34136 Trieste, Italy}
} 
\bigskip

\centerline{\scshape Alessio Figalli }
\medskip
{\footnotesize
 \centerline{University of Texas at Austin}
 \centerline{Department of Mathematics}
   \centerline{2515 Speedway Stop C1200}
   \centerline{Austin, TX 78712-1202, USA}
}
\bigskip

\centerline{\scshape Giampiero Palatucci }
\medskip
{\footnotesize
\centerline{Universit\`a degli Studi di Parma}
 \centerline{Dipartimento di Matematica}
\centerline{Campus - Parco Area delle Scienze, 53/A}
\centerline{43124 Parma, Italy}
}
\bigskip

\centerline{\scshape Enrico Valdinoci }
\medskip
{\footnotesize
\centerline{Universit\`a degli Studi di Milano}
 \centerline{Dipartimento di Matematica}
\centerline{Via Saldini, 50}
\centerline{20133 Milano, Italy}
}
\bigskip

\centerline{(Communicated by the associate editor name)}

\hyphenation{par-ma ga-gli-ar-do}

\begin{abstract}
We deal with the asymptotic behavior of the $s$-perimeter of a set~$E$ inside a domain $\Omega$
as $s\searrow0$.
We prove necessary and sufficient conditions for the existence of such limit,  
by also providing an explicit formulation in terms of the Lebesgue measure of $E$ and $\Omega$. 
Moreover, we construct examples of sets for which the limit does not exist.
\end{abstract}

\section{Introduction}

Given~$s\in(0,1)$ and a bounded open set $\Omega\subset\R^n$
with $C^{1,\gamma}$-boundary,
the $s$-perimeter
of a (measurable) set $E\subseteq\R^n$ in $\Omega$ 
is defined as
\begin{equation}\label{def_per}
\begin{split}
&{\text{\rm Per}}_{s} (E;\Omega):=L(E\cap\Omega, (\CC E)\cap\Omega)
\\ &\qquad\qquad\qquad+
L(E\cap\Omega, (\CC E)\cap(\CC\Omega))+
L(E\cap(\CC\Omega), (\CC E)\cap\Omega),
\end{split}
\end{equation}
where $\CC E = \R^n\setminus E$ denotes the complement of~$E$, 
and $L(A,B)$ denotes the following nonlocal interaction term
\begin{equation}\label{def_l}
\displaystyle
L(A,B):=\int_A \int_B \frac{1}{|x-y|^{n+s}}\,dx\,dy \qquad \forall\, A, B\subseteq\R^n.
\end{equation}
Here we are using the standard convention
for which $L(A,B)=0$ if either $A=\varnothing$ or $B=\varnothing$.

This notion of~$s$-perimeter and the corresponding 
minimization problem were introduced in~\cite{CRS10}
(see also the pioneering work~\cite{Vis87, Vis91},
where some functionals related
to the one in~\eqref{def_per} have been analyzed 
in connection with fractal dimensions).

Recently, the $s$-perimeter has inspired
a variety of literature in different directions, 
both in the pure mathematical settings (for instance,
as regards the regularity of surfaces with minimal
$s$-perimeter, see~\cite{BFV12, CG12,CV12, SV12}) 
and in view of concrete applications (such as
phase transition 
problems with long range interactions, see~\cite{CS10, SV11, SV11b}).
In general, the nonlocal behavior of the functional is the
source of major difficulties, conceptual differences, and challenging technical
complications. We refer to~\cite{FV12} for an introductory review on this subject.

The limits as~$s\searrow0$ and~$s\nearrow1$ are somehow the critical cases for the 
$s$-perimeter, since the functional in~\eqref{def_per} diverges as it is. 
Nevertheless, when appropriately rescaled, these limits seem to give meaningful 
information on the problem. In particular, 
it was shown in~\cite{CV11,ADM11} that $(1-s)
{\text{\rm Per}}_{s}$
approaches the classical perimeter functional as~$s\nearrow1$  (up
to normalizing multiplicative constants),
and this implies that surfaces of minimal $s$-perimeter
inherit the regularity properties of the classical minimal surfaces
for~$s$ sufficiently close to~$1$ (see~\cite{CV12}).

As far as we know, the asymptotic as~$s\searrow0$
of~$s{\text{\rm Per}}_{s}$ was not studied yet (see however~\cite{MS}
for some results in this direction), and this is the question
that we would like to address in this paper.
That is, we are interested in the quantity
\begin{equation}\label{RR}
 \mu(E):=\lim_{s\searrow0} 
s{\text{\rm Per}}_{s}(E;\Omega)
\end{equation}
whenever the limit exists. 
Of course,
if it exists then
$$ \mu(E)=\mu(\CC E),$$ 
since $$
{\text{\rm Per}}_{s}(E;\Omega)=
{\text{\rm Per}}_{s}(\CC E;\Omega).$$

We will show that,
though~$\mu$ is subadditive (see Proposition~\ref{T1} below), in
general it is not a 
measure (see Proposition~\ref{T2}, and this is a major difference with
respect to the setting in~\cite{MS}). 
On the other hand,~$\mu$ is additive on bounded, separated sets, and 
it agrees with the Lebesgue measure of $E \cap \Omega$ (up to normalization) when~$E$ is bounded
(see Corollary~\ref{Cor}). 
As we will show below, a precise characterization of~$\mu(E)$ will be given in terms of 
the behavior of the set~$E$ towards infinity, which is encoded
in the quantity
$$ \alpha(E):=\lim_{s\searrow0} s\int_{E\cap (\CC B_1)}\frac{1}{|y|^{n+s}}\,dy,$$
whenever it exists (see Theorem~\ref{TF}
and Corollary~\ref{Cor}).
In fact, the existence of the limit defining~$\alpha$
is in general equivalent to the one defining~$\mu$ (see 
Theorem~\ref{TF1}(ii)).

As a counterpart of these results, we will construct an explicit example of set $E$ for which both the 
limits~$\mu(E)$ and~$\alpha(E)$ do not exist (see 
Example~\ref{EXX}): this says that the assumptions we take cannot, in general,
be removed.

Also, notice that, in order to make sense of the limit in~\eqref{RR}, it is
necessary to assume that\footnote{It is easily seen that if~\eqref{s0} holds,
then~$ {\text{\rm Per}}_{s}(E;\Omega)<\infty$ for any~$s\in(0,s_0)$.
Moreover, if~$\partial E$ is smooth, then~\eqref{s0} is always satisfied.}
\begin{equation}\label{s0}
{\mbox{$ {\text{\rm Per}}_{s_0}(E;\Omega)<\infty$, for some $s_0 \in (0,1)$.}}
\end{equation}
To stress that~\eqref{s0} cannot be dropped,
we will construct a simple example in which such a condition is violated
(see Example~\ref{EX2}).
\vspace{3mm}

The paper is organized as follows. In the following section,
we collect the precise statements of all the results we mentioned above.
Section~\ref{sec_proofs} is devoted to the proofs.

\vspace{2mm}

\section{List of the main results} 

We define~${\mathcal{E}}$ to be the family of sets~$E\subseteq\R^n$ for which
the limit defining~$\mu(E)$ in~\eqref{RR} exists.
We prove the following result:

\begin{prop}\label{T1}
$\mu$ is subadditive on~${\mathcal{E}}$, i.e. $\mu(E\cup F)\le \mu(E)+\mu(F)$
for any $E$, $F\in{\mathcal{E}}$.
\end{prop}

First, it is convenient to consider the normalized
Lebesgue measure~${\mathcal{M}}$, that is
the standard Lebesgue measure scaled by the factor
${\mathcal{H}}^{n-1}(S^{n-1})$, namely
\begin{equation}\label{2.1bis} {\mathcal{M}}(E):=
{\mathcal{H}}^{n-1}(S^{n-1})\,|E|,
\end{equation}
where, as usual, we denote by $S^{n-1}$ the $(n\!-\!1)$-dimensional sphere.
\vspace{2mm}

Now, we recall the main result in~\cite{MS}; that is,
\begin{theorem}\label{thm_ms}{\rm (see \cite[Theorem 3]{MS}).}
Let $s\in(0,1)$. Then, for all $u\in H^s(\R^n)$,
\begin{equation*}
\displaystyle
\lim_{s \searrow0} \,\frac{s}{2}\,\int_{\R^n}\int_{\R^n} \frac{|u(x)-u(y)|^2}{|x-y|^{n+s}}\,dx\,dy
\, = \, \mathcal{H}^{n-1}({S}^{n-1}) \int_{\R^n} |u|^2\, dx.
\end{equation*}
\end{theorem}
An easy consequence of the result above is that when~$E\in {\mathcal{E}}$ and~$E\subseteq\Omega$
then~$\mu(E)$ agrees with~${\mathcal{M}}(E)$
(in fact, we will generalize this statement in Theorem~\ref{TF}
and Corollary~\ref{Cor}). Based on this property valid for subsets
of~$\Omega$, one may be tempted to infer that~$\mu$ is always
related to the Lebesgue measure, up to normalization,
or at least to some more general type of measures.
The next result points out that this cannot be true:

\begin{prop}\label{T2}
$\mu$ is not necessarily additive on separated sets in~$ {\mathcal{E}}$, i.e. there exist
$E,F\in {\mathcal{E}}$ such that 
$\text{\rm dist}(E,F)\ge c>0$, but
$\mu(E\cup F)< \mu(E)+\mu(F)$. 

Also, $\mu$ is not necessarily
monotone on~${\mathcal{E}}$, i.e. it is not true that $E\subseteq F$ implies $\mu(E)\le\mu(F)$.
\end{prop}

In particular, we deduce
from Proposition~\ref{T2}
that $\mu$ is not a measure.
On the other hand, in some circumstances
the additivity property holds true:

\begin{prop}\label{T3}
$\mu$ is additive on bounded, separated sets
in~${\mathcal{E}}$, i.e. if~$E$, $F\in{\mathcal{E}}$, $E$ and $F$ are bounded, disjoint
and \,$\text{\rm dist}(E,F)\ge c>0$, then~$E\cup F\in{\mathcal{E}}$ and~$\mu(E\cup F)=\mu(E)+\mu(F)$.
\end{prop}

There is a natural condition under which~$\mu(E)$ does 
exist, based on the weighted volume of~$E$ towards infinity,
as next result points out:

\begin{theorem}\label{TF}
Suppose that \,${\text{\rm Per}}_{s_0}(E;\Omega)<\infty$ for some $s_0 \in (0,1)$, 
and that the 
following 
limit exists
\begin{equation}\label{Rj} 
\alpha(E):=\lim_{s\searrow0}
s\int_{E\cap (\CC B_1)}\frac{1}{|y|^{n+s}}\,dy.\end{equation}
Then~$E\in{\mathcal{E}}$ and
$$ 
\mu(E)=\big(1-\widetilde\alpha(E)\big)\,{\mathcal{M}}(E\cap\Omega)+
\widetilde\alpha(E)\,{\mathcal{M}}(\Omega\setminus E),$$\
where
\begin{equation}\label{07}
\widetilde\alpha(E):=\frac{\alpha(E)}{
{\mathcal{H}}^{n-1}(S^{n-1})}.\end{equation}
\end{theorem}  

\vspace{2mm}

As a consequence of Theorem \ref{TF}, one obtains the existence and the exact 
expression of $\mu(E)$ for a bounded set $E$, as described by the following result:

\begin{corollary}\label{Cor}
Let~$E$ be a bounded set, 
and \,${\text{\rm Per}}_{s_0}(E;\Omega)<\infty$ for some $s_0 \in (0,1)$. 
Then~$E\in{\mathcal{E}}$ and
$$ 
\mu(E)={\mathcal{M}}(E\cap\Omega).$$
In particular,
if~$E\subseteq\Omega$ and~${\text{\rm Per}}_{s_0}(E;\Omega)<\infty$ for some $s_0 \in (0,1)$,
then $\mu(E)={\mathcal{M}}(E)$.
\end{corollary}

Condition \eqref{Rj} is also in general necessary for the existence of
the limit in~\eqref{RR}. Indeed, next result shows that the existence of
the limit in \eqref{Rj} is equivalent to the existence of the limit in \eqref{RR},
except in the special case in which the set $E$ occupies exactly half of the measure
of $\Omega$ (in this case the limit in
\eqref{RR} always exists, independently on the existence of the limit in \eqref{Rj}).

\begin{theorem}\label{TF1}
Suppose that \,${\text{\rm Per}}_{s_0}(E;\Omega)<\infty$, for some $s_0 \in (0,1)$.
Then:
\begin{enumerate}
\item[{(i)}] If $|\Omega\setminus  E| = |E \cap \Omega|$, then $E\in{\mathcal{E}}$ and $\mu(E)=\mathcal M(E\cap \Omega)$.
\item[{(ii)}] If $|\Omega\setminus  E| \neq |E \cap \Omega|$
and~$E\in{\mathcal{E}}$, then the 
limit in~\eqref{Rj} exists 
and 
$$ \alpha(E) = \frac{\mu(E) - {\mathcal{M}}(E \cap \Omega)}{
|\Omega\setminus  E| - |E \cap \Omega|}. $$
\end{enumerate}
\end{theorem}  
\vspace{2mm}

In the statements above we assumed the existence
of the limits in~\eqref{RR} and~\eqref{Rj}. Such assumptions
cannot be removed, since
the limits in \eqref{RR} and~\eqref{Rj} may not exist, as we
now point out:

\begin{example}\label{EXX} There exists a set $E$ with $C^\infty$-boundary for which the limits
in \eqref{RR} and 
\eqref{Rj} do not exist.
\end{example}

\begin{example}\label{EXX special} There exists a set $E$ with
$C^\infty$-boundary for which the limit
in \eqref{RR} exists and the limit in
\eqref{Rj} does not exist.
\end{example}

Notice that Examples \ref{EXX} and \ref{EXX special}
are provided by smooth sets, and therefore they have
finite $s$-perimeter for any $s\in(0,1)$ (see, e.g.,
Lemma 11 in \cite{CV11}).

On the other hand, as regards condition~\eqref{s0}, we point out
that it cannot be dropped in general, since there are
sets that do not satisfy it (and for them the limit
in~\eqref{RR} does not make sense):

\begin{example}\label{EX2} There exists a set $E$ 
for which \,$\text{\rm Per}_s(E;\Omega)=+\infty$ for any $s\in(0,1)$.
\end{example}

\vspace{2mm}

\section{Proofs}\label{sec_proofs} 

\subsection{Proof of Proposition \ref{T1}}\label{PF T1}

We observe that
\begin{equation}\label{sub s}
{\mbox{the $s$-perimeter is subadditive.}}
\end{equation}
To check this, let $\Omega_1$, $\Omega_2$ be open sets
of $\R^n$. We remark that
\begin{eqnarray*}
&& \!\!\!\!\!\!L((E\cup F)\cap\Omega_1,(\CC (E\cup F))\cap \Omega_2)
\\[1ex]
&&\qquad \quad=
L((E\cap\Omega_1)\cup(F\cap\Omega_1),(\CC E)\cap(\CC F)\cap\Omega_2)\\[1ex]
&&\qquad\quad\le
L( E\cap\Omega_1,(\CC E)\cap(\CC F)\cap \Omega_2)
+
L(F\cap\Omega_1,(\CC E)\cap(\CC F)\cap \Omega_2)
\\[1ex]
&&\qquad\quad\le
L( E\cap\Omega_1,(\CC E)\cap \Omega_2)
+
L(F\cap\Omega_1,(\CC F)\cap \Omega_2).
\end{eqnarray*}
By taking $\Omega_1:=\Omega$ and $\Omega_2:=\R^n$
we obtain
$$ L((E\cup F)\cap\Omega,\CC (E\cup F))\le 
L( E\cap\Omega,\CC E)
+
L(F\cap\Omega,\CC F),$$
while, by taking $\Omega_1:=\CC\Omega$ and $\Omega_2:=\Omega$,
we conclude that
$$ L((E\cup F)\cap(\CC\Omega),(\CC (E\cup F))\cap \Omega)\leq
L( E\cap(\CC\Omega),(\CC E)\cap \Omega)
+
L(F\cap(\CC\Omega),(\CC F)\cap \Omega).$$
By summing up, we get
\begin{eqnarray*}
&& \!\!\!\!\!\!\!\!\!\!\!\!
{\text{\rm Per}}_s(E\cup F;\Omega)\\
&& \quad = L( (E\cup F)\cap\Omega, \CC (E\cup F))+
L((E\cup F)\cap(\CC\Omega),(\CC(E\cup F))\cap\Omega)
\\[1ex]
&& \quad\le
L( E\cap\Omega,\CC E)
+
L(F\cap\Omega,\CC F)\\
&&\quad \quad +\,
L( E\cap(\CC\Omega),(\CC E)\cap \Omega)
+
L(F\cap(\CC\Omega),(\CC F)\cap \Omega)
\\[1ex]
&& \quad =
\text{\rm Per}_s(E;\Omega)+ \text{\rm Per}_s(F;\Omega)
.\end{eqnarray*}
This establishes \eqref{sub s} and then Proposition \ref{T1} follows by taking the limit as $s\searrow0$.~\hfill$\square$

\subsection{Proof of Proposition \ref{T2}}\label{PF T2}

First we show that~$\mu$ is not additive.

Here and in the sequel, we denote by $B_R$ the open ball centered at $0\in \R^n$ of radius $R>0$. 
We observe that if $x\in B_1$ and $y\in\CC B_2$ then $|x-y|\,\le\,|x|+|y|\,\le\, 2|y|$, therefore
$$ s L(B_1,\CC B_2)\,\ge\, c_1 s
\int_{B_1} dx \int_{\CC B_2} dy \frac{1}{|y|^{n+s}}\,\ge\, c_2 s\int_2^{+\infty}
\frac{d\rho}{\rho^{1+s}}
\,\ge\, c_3,$$
for some positive constants $c_1$, $c_2$ and $c_3$.
Now we take $E:=\CC B_2$, $F:=\Omega:=B_1$. Then
\begin{eqnarray*}
&& {\text{\rm Per}}_{s} (E;\Omega)=L(B_1,\CC B_2),\\
&& {\text{\rm Per}}_{s} (F;\Omega)=L(B_1,\CC B_1)=L(B_1,\CC B_2)+L(B_1,B_2\setminus B_1)\\
\aand && {\text{\rm Per}}_{s} (E\cup F;\Omega)=
L(B_1,B_2\setminus B_1).
\end{eqnarray*}
Therefore
\begin{eqnarray*}
s\,{\text{\rm Per}}_{s} (E;\Omega)
+s\,{\text{\rm Per}}_{s} (F;\Omega) &= & 2sL(B_1,\CC B_2)+sL(B_1,B_2\setminus B_1) \\
&  \ge &2c_3+s\,
L(B_1,B_2\setminus B_1)
\\
& = & 2c_3 +s\,{\text{\rm Per}}_{s} (E \cup F;\Omega).
\end{eqnarray*}
By sending $s\searrow0$, we conclude that $\mu(E)+\mu(F)\ge 2c_3+\mu(E\cup F)$, so~$\mu$ is not additive.

Now we show that~$\mu$ is not monotone either.
For this we take~$E$ 
such that $\mu(E)>0$ (for instance,
one can take $E$ a small ball inside $\Omega$; see Corollary~\ref{Cor}),
and~$F:=\R^n$:
with this choice,~$E\subset F$ and~$\text{\rm Per}_s(F;\Omega)=0$,
so~$\mu(E)>0=\mu(F)$.~\hfill$\square$
\vspace{2mm}

\subsection{Auxiliary observations}

Here we collect some observations, to be exploited in the 
subsequent proofs. 
\vspace{1mm}

\noindent{\bf Observation 1.}
First of all, 
we observe that
\begin{equation}\label{T3 E1}\begin{split}
&{\mbox{if $A$ and $B$ are bounded, disjoints sets with~$\dist(A,B)\ge 
c>0$, then }}\\
&\qquad \displaystyle\lim_{s\searrow0} s \,L(A,B)=0.
\end{split}\end{equation}
To check this, suppose that $A$ and $B$ lie in $B_R$. Then
$$ \int_A \int_B \frac{1}{|x-y|^{n+s}}\,dx\,dy
\,\le\, 
\int_{B_R} \int_{B_R} \frac{1}{c^{n+s}}\,dx\,dy
\,=\, 
\frac{(\mathcal{H}^{n-1}({S}^{n-1}))^2 R^{2n} }{n^2c^{n+s}}$$
and this establishes \eqref{T3 E1}.
\smallskip

\noindent{\bf Observation 2.}
Now we would like to
remark that the quantity
$$ \lim_{s\searrow0}
s\int_{E\cap (\CC B_R)}\frac{1}{|y|^{n+s}}\,dy$$
is independent of~$R$, if the limit exists.
More precisely, we show that
for any~$R\ge r>0$ 
\begin{equation}\label{0.01}
\lim_{s\searrow0}
s\left(\int_{E\cap (\CC B_R)}\frac{1}{|y|^{n+s}}\,dy
-\int_{E\cap (\CC B_r)}\frac{1}{|y|^{n+s}}\,dy\right)=0.
\end{equation}
To prove this, we notice that
\begin{eqnarray}
   s\int_{E \cap (B_R \setminus  B_r)}\frac{1}{|y|^{n+s}}\,dy 
   \leq     s\int_{B_R \setminus  B_r}\frac{1}{|y|^{n+s}}\,dy  
   = s {\mathcal{H}}^{n-1}(S^{n-1}) \int_{r}^{R} 
\frac{1}{\rho^{1+s}}\,d\rho \nonumber\\
   = {\mathcal{H}}^{n-1}(S^{n-1})\left( \frac{1}{r^s} 
-\frac{1}{R^s}\right) 
\end{eqnarray}
and so, by taking limit in $s$,
\begin{equation*}
\lim_{s\searrow0} s\int_{E \cap (B_R \setminus  B_r)}
\frac{1}{|y|^{n+s}}\,dy  =0, 
\end{equation*}
which establishes~\eqref{0.01}.\smallskip

\noindent{\bf Observation 3.}
As a consequence of~\eqref{0.01}, it follows that if the limit
in~\eqref{Rj} exists then
\begin{equation}\alpha(E)= \lim_{s\searrow0}
s\int_{E\cap (\CC B_R)}\frac{1}{|y|^{n+s}}\,dy \qquad \forall\,R>0. \label{alfa}\end{equation}
\smallskip

\noindent{\bf Observation 4.}
For any $s\in(0,1)$, we define
\begin{equation}\label{AS}
\alpha_s(E) := s \int_{E \cap (\CC B_1)}  \frac{1}{|y|^{n+s}}\,dy \end{equation}
and we prove that, for any bounded set~$F\subset\R^n$, and any set~$E\subseteq\R^n$,
\begin{equation}\label{XX00}
\lim_{R\rightarrow+\infty}
\limsup_{s\searrow0}\left|\alpha_s(E)\,|F| -s \int_{F}\int_{E\cap 
(\CC B_R)}\frac{1}{|x-y|^{n+s}}\,dx\,dy\right|=0.\end{equation}
To prove this, we take~$r>0$ such that~$F\subset B_r$ and~$R>1+2r$
(later on~$R$ will be taken as large as we wish). 
We observe that, for any~$z\in B_r$ and~$y\in\CC B_R$,
$$ |z-y|\ge |y|-|z|=\left(1-\frac{r}{R}\right)|y|+\frac{r}{R}|y|-|z|\ge
\frac{|y|}2.$$
Therefore, if, for any fixed~$y\in\CC B_R$ we consider the map
$$ h(z):=\frac{1}{|z-y|^{n+s}},\qquad z\in B_r,$$
we have that
$$ |\nabla h(z)|=\frac{n+s}{|z-y|^{n+s+1}}\le \frac{2^{n+s+1}(n+s)}{|y|^{n+s+1}},$$
for any~$z\in B_r$, which implies 
$$ \left| \frac{1}{|x-y|^{n+s}}-\frac{1}{|y|^{n+s}}\right|=|h(x)-h(0)|\le
\frac{2^{n+s+1}(n+s)|x|}{|y|^{n+s+1}}\qquad \forall\,x\in B_r,\,y\in\CC B_R.
$$
Therefore
\begin{eqnarray*}&&
\left|\int_F\left(\int_{E\cap(\CC B_R)}\frac{1}{|y|^{n+s}}\,dy\right)\,dx-\int_{F}\int_{E\cap
(\CC B_R)}\frac{1}{|x-y|^{n+s}}\,dx\,dy\right|\\
&&\qquad \le 
\int_F\left(\int_{E\cap(\CC B_R)}\left|\frac{1}{|y|^{n+s}}
-\frac{1}{|x-y|^{n+s}} \right|\,dy\right)\,dx\\
&&\qquad\le
\int_F\left(
\int_{E\cap(\CC B_R)}
\frac{2^{n+s+1}(n+s)|x|}{|y|^{n+s+1}}
\,dy\right)\,dx
\\&&\qquad\le
2^{n+s+1}(n+s) |F| r
\int_{\CC B_R}
\frac{1}{|y|^{n+s+1}}\,dy\,\le\,C
\end{eqnarray*}
for some~$C>0$ independent of $s$. As a consequence
\begin{eqnarray*}
&& \left|\alpha_s(E)\,|F| -s \int_{F}\int_{E\cap
(\CC B_R)}\frac{1}{|x-y|^{n+s}}\,dx\,dy\right|\\
&&\qquad\qquad\qquad\qquad\qquad\le |F|\,\left|\alpha_s(E)-s\int_{E\cap(\CC B_R)}
\frac{1}{|y|^{n+s}}\,dy\right|+Cs.
\end{eqnarray*}
This and~\eqref{0.01} (applied here with~$r:=1$)
imply~\eqref{XX00}.
\smallskip

\noindent{\bf Observation 5.} If the limit in~\eqref{Rj}
exists, then~\eqref{XX00} boils down to
\begin{equation}\label{79}
\lim_{R\rightarrow+\infty}
\limsup_{s\searrow0}\left|\alpha(E)\,|F| -s \int_{F}\int_{E\cap
(\CC B_R)}\frac{1}{|x-y|^{n+s}}\,dx\,dy\right|=0.
\end{equation}
\smallskip

\noindent{\bf Observation 6.}
Now we point out
that, if 
$F\subseteq\Omega\subset B_R$ for some $R>0$, and~$F$ has finite $s_0$-perimeter
in~$\Omega$ for some $s_0 \in (0,1)$, 
then
\begin{equation}\label{700}
\lim_{s\searrow0} s\int_{F}\int_{B_R\setminus
F}\frac{1}{|x-y|^{n+s}}\,dx\,dy=0.
\end{equation}
Indeed, for any~$s\in(0,s_0)$,
\begin{eqnarray*}
&& \int_{F}\int_{B_R\setminus
F}\frac{1}{|x-y|^{n+s}}\,dx\,dy\\ &&\qquad\le
\int_{F}\int_{(B_R\setminus F)\cap\{|x-y|\le 1\} }
\frac{1}{|x-y|^{n+s_0}}\,dx\,dy
+
\int_{F}\int_{(B_R\setminus F)\cap\{|x-y|> 1\}}
1\,dx\,dy
\\ &&\qquad\le {\text{\rm Per}}_{s_0}(F;\Omega)+|B_R|^2,\end{eqnarray*}
which implies~\eqref{700}. In particular, 
thanks to \cite[Proposition 16]{ADM11}, the argument above also shows that
if $F\Subset\Omega\subset B_R$ 
and $\chi_F \in BV(\Omega)$, then $F$ has finite $s$-perimeter
in~$\Omega$ for any $s \in (0,1)$.
\smallskip

\noindent{\bf Observation 7.}
Let~$E_1:=E\cap
\Omega$ and~$E_2:=E\setminus \Omega$. Then
\begin{equation}\label{contribution}\begin{split}
{\text{\rm Per}}_s(E;\Omega)\,&=
{\text{\rm Per}}_s(E_1\cup E_2;\Omega)
\\ &= L(E_1,\Omega\setminus E_1)+
L(E_1,(\CC \Omega)\setminus E_2)+L(E_2,\Omega\setminus E_1)
\\ &= L(E_1,\CC E_1)-L(E_1,E_2)+L(E_2,\Omega\setminus E_1)
\\ &= {\text{\rm Per}}_s(E_1;\Omega)-L(E_1,E_2)+L(E_2,\Omega\setminus E_1)
.\end{split}\end{equation}

\medskip

With these observations in hand, we are ready to continue the proofs
of the main results.

\subsection{Proof of Proposition \ref{T3}}\label{PF T3}

We prove Proposition \ref{T3}
by suitably modifying the proof of Proposition \ref{T1}.
Given two open sets
$\Omega_1$ and $\Omega_2$, and two disjoint sets $E$ and $F$,
we have that
\begin{eqnarray*}
&&\!\!\!\!\!\!\!\! L((E\cup F)\cap\Omega_1,(\CC (E\cup F))\cap \Omega_2)
\\[1ex]
&&\qquad=
L((E\cap\Omega_1)\cup(F\cap\Omega_1),(\CC E)\cap(\CC F)\cap\Omega_2)\\[1ex]
&&\qquad =
L(E\cap\Omega_1,(\CC E)\cap(\CC F)\cap\Omega_2)
+
L(F\cap\Omega_1,(\CC E)\cap(\CC F)\cap\Omega_2).
\end{eqnarray*}
By taking $\Omega_1:=\Omega$ and $\Omega_2:=\R^n$
we obtain
$$ L((E\cup F)\cap\Omega, \CC (E\cup F) )=
L(E\cap\Omega,(\CC E)\cap(\CC F))
+ L(F\cap\Omega,(\CC E)\cap(\CC F) )$$
while, by taking $\Omega_1:=\CC\Omega$ and $\Omega_2:=\Omega$,
we conclude that
\begin{eqnarray*}  
&& \!\!\!\!\!\!\!\! L((E\cup F)\cap(\CC\Omega),(\CC (E\cup F))\cap \Omega) \\
&& \qquad = \, L(E\cap(\CC\Omega),(\CC E)\cap(\CC F)\cap\Omega) 
 +
L(F\cap(\CC\Omega),(\CC E)\cap(\CC F)\cap\Omega).\end{eqnarray*}
As a consequence,
\begin{eqnarray*}
&& \!\!\!\!\!\!
{\text{\rm Per}}_s(E\cup F;\Omega)
\\ &&\qquad = \, L((E\cup F)\cap\Omega, \CC (E\cup F) )+
L((E\cup F)\cap(\CC\Omega),(\CC (E\cup F))\cap \Omega)\\[1ex]
 &&\qquad = \, L(E\cap\Omega,(\CC E)\cap(\CC F))
+ L(F\cap\Omega,(\CC E)\cap(\CC F) )\\
 &&\qquad \quad + \,
L(E\cap(\CC\Omega),(\CC E)\cap(\CC F)\cap\Omega)
+
L(F\cap(\CC\Omega),(\CC E)\cap(\CC F)\cap\Omega)
\\[1ex]
 &&\qquad = \,{\text{\rm Per}}_s(E;\Omega)+
{\text{\rm Per}}_s(F;\Omega)
\\ &&\qquad \quad -\, L(E\cap\Omega, (\CC E)\cap F)
-L(F\cap\Omega, E\cap(\CC F))
\\ &&\qquad \quad -\,L(E\cap(\CC\Omega), (\CC E)\cap F\cap\Omega)
-L(F\cap(\CC\Omega), E\cap(\CC F)\cap\Omega).
\end{eqnarray*}
We remark that the last interactions involve only 
bounded, separated sets, since so are~$E$ and~$F$, therefore, by \eqref{T3 E1},
$$ \lim_{s\searrow0} s \,{\text{\rm Per}}_s(E\cup F;\Omega)=
\lim_{s\searrow0} \big( s \,{\text{\rm Per}}_s(E;\Omega)
+s \,
{\text{\rm Per}}_s(F;\Omega)\big),$$ 
which completes the proof of Proposition \ref{T3}.~\hfill$\square$

\subsection{Proof of Theorem \ref{TF}}\label{PF TF}

We suppose that~$\Omega\subset B_r$, for some~$r>0$,
and we take ~$R>1+2r$. 
Let~$E_1:=E\cap 
\Omega$ and~$E_2:=E\setminus \Omega$.
Notice that, for any~$F\subseteq\Omega$, which has finite $s_0$-perimeter in $\Omega$ for some $s_0\in(0,1)$, 
$$E_2\cap B_R\subseteq B_R\setminus\Omega\subseteq
B_R\setminus F$$
and so~\eqref{700} gives that
\begin{equation} \label{7.-1}
\lim_{s\searrow0} s\int_{F}\int_{E_2\cap 
B_R}\frac{1}{|x-y|^{n+s}}\,dx\,dy=0,\end{equation}
provided that $F$ has finite $s_0$-perimeter in $\Omega$.
Using this and~\eqref{79}, we conclude that, for any~$F\subseteq\Omega$
of finite $s_0$-perimeter in $\Omega$,
\begin{equation*}\begin{split}
&\lim_{s\searrow0}
s\int_{F}\int_{E_2}\frac{1}{|x-y|^{n+s}}\,dx\,dy
\\&\qquad\qquad\qquad\qquad=\,
\lim_{R\rightarrow+\infty}\lim_{s\searrow0}
s\int_{F}\int_{E_2}\frac{1}{|x-y|^{n+s}}\,dx\,dy
\\ &\qquad\qquad\qquad\qquad=\,\lim_{R\rightarrow+\infty}\lim_{s\searrow0}
s\int_{F}\int_{E_2\cap(\CC B_R)}\frac{1}{|x-y|^{n+s}}\,dx\,dy
\\[1ex]
&\qquad\qquad\qquad\qquad=\,\alpha(E)\,|F|.\end{split}
\end{equation*}
In particular\footnote{We stress that both $E_1$ and $\Omega\setminus E_1$
have finite $s_0$-perimeter in $\Omega$ if so has $E$, thanks
\label{foot}
to our smoothness assumption on $\partial\Omega$.
We check this claim for $E_1$, the other being analogous. First
of all, fixed $B_R\supset B_r\supset\Omega$, we have that
$$ L\big(E_1,(E\setminus\Omega)\cap(\CC B_R)\big)\le
L(B_r, \CC B_R)<+\infty.$$
Also $L\big(\Omega \cap B_R, (\CC\Omega)\cap B_R\big)<+\infty$
(see, e.g., Lemma 11 in \cite{CV11}), therefore
\begin{eqnarray*} && {\text{\rm Per}}_{s_0} (E_1;\Omega)=L(E_1,\CC E_1)=
L(E_1,\CC E)+L(E_1, E\setminus\Omega)\\ &&\qquad\le
{\text{\rm Per}}_{s_0} (E;\Omega)
+L\big(E_1, (E\setminus\Omega)\cap B_R\big)+
L\big(E_1, (E\setminus\Omega)\cap (\CC B_R)\big)
\\ &&\qquad\le
{\text{\rm Per}}_{s_0} (E;\Omega)
+L\big(\Omega, (\CC\Omega)\cap B_R\big)+
L\big(E_1, (E\setminus\Omega)\cap (\CC B_R)\big),
\end{eqnarray*}
that is finite.}, by taking~$F:=E_1$ and~$F:=\Omega\setminus E_1$, 
and recalling \eqref{2.1bis} and~\eqref{07},
\begin{equation}\label{90}
\begin{split}
&\lim_{s\searrow0}
s\int_{E_1}\int_{E_2}\frac{1}{|x-y|^{n+s}}\,dx\,dy=\alpha(E)\,|E_1|
=\widetilde\alpha(E)\,{\mathcal{M}}(E_1)
\\
{\mbox{and }}\;&
\lim_{s\searrow0}
s\int_{\Omega\setminus E_1}\int_{E_2}\frac{1}{|x-y|^{n+s}}\,dx\,dy=
\alpha(E)\,|\Omega\setminus 
E_1|=\widetilde\alpha(E)\,{\mathcal{M}}(\Omega\setminus E_1)
.\end{split}\end{equation}
We now claim 
\begin{equation}\label{E1}\lim_{s\searrow0}
s\,{\text{\rm Per}}_s(E_1;\Omega)={\mathcal{M}}(E_1).
\end{equation}
Indeed, since $E_1\subseteq \Omega$, this
is a plain consequence of Theorem~\ref{thm_ms} (see
also Remark~4.3 in \cite{DPV12} for another elementary proof) by simply choosing $u=\chi_{E_1}$ there:
\begin{equation*}\begin{split}
\lim_{s\searrow0} s \,
{\text{\rm Per}}_s(E_1;\Omega)
& = \, \lim_{s\searrow0} s \,L(E_1, \CC E_1) \\
& = \, \lim_{s\searrow0} \frac{s}{2}\int_{\R^n}\int_{\R^n}\frac{ |\chi_{E_1}(x)-\chi_{E_1}(y)|^2}{|x-y|^{n+s}}\,dx\,dy
\\[1ex]
& = \,{\mathcal{H}}^{n-1}(S^{n-1})\,\| \chi_{E_1}\|_{L^2(\R^n)}^2 \ = \ 
{\mathcal{H}}^{n-1}(S^{n-1})\,|E_1|,
\end{split}\end{equation*}
as desired.
Thus, using~\eqref{contribution}, \eqref{90}, and~\eqref{E1}, we 
obtain
\begin{eqnarray*}
\lim_{s\searrow0} s
{\text{\rm Per}}_s(E;\Omega)
= {\mathcal{M}}(E_1)-\widetilde\alpha(E) {\mathcal{M}}(E_1)+
\widetilde\alpha(E) {\mathcal{M}}(\Omega\setminus E_1),
\end{eqnarray*}
which is the desired result.~\hfill$\square$

\subsection{Proof of Corollary \ref{Cor}}\label{PF Cor}
We fix $R$ large enough so that $E\subset B_R$, 
hence $E \cap (\CC B_R) = \varnothing$. 
By the expression of $\alpha(E)$ in $(\ref{alfa})$, 
we have that the limit in~$(\ref{Rj})$ exists and $\alpha(E)=0$. 
Then the result follows by Theorem \ref{TF}.~\hfill$\square$

\subsection{Proof of Theorem \ref{TF1}}
We suppose that~$\Omega\subset B_r$, for some~$r>0$, 
and we take ~$R>1+2r$.
Let~$E_1:=E\cap\Omega$ 
and~$E_2:=E\setminus \Omega$.
By~\eqref{contribution},
\begin{equation*}
\begin{split} 
&\!\!\!\!\! s 
{\text{\rm Per}}_s(E;\Omega) -s {\text{\rm Per}}_s(E_1;\Omega) \\[1ex]
& \qquad=
s L(E_2,\Omega\setminus E_1) - s L(E_1,E_2) \\[1ex]
 & \qquad = s \int_{\Omega\setminus E_1} \int_{E_2 \cap B_R} 
\frac{1}{|x-y|^{n+s}}\, dx \, dy  
  +  s \int_{\Omega\setminus E_1} \int_{E_2 \cap (\CC B_R)} 
\frac{1}{|x-y|^{n+s}}\, dx \, dy  \\
& \qquad 
 \quad-   s \int_{ E_1} \int_{E_2 \cap B_R} \frac{1}{|x-y|^{n+s}}\, dx \, dy    
 -   s \int_{ E_1} \int_{E_2 \cap (\CC B_R)} \frac{1}{|x-y|^{n+s}}\, dx 
\, dy. 
\end{split}\end{equation*}
By rearranging the terms, we obtain
\begin{equation}\begin{split}\label{10}
& I(s,R):=s \int_{\Omega\setminus E_1} \int_{E_2 \cap (\CC B_R)} 
\frac{1}{|x-y|^{n+s}}\, dx \, dy  
- s \int_{ E_1} \int_{E_2 \cap (\CC B_R)} \frac{1}{|x-y|^{n+s}}\, dx \, 
dy\\[1ex]
& \qquad \quad \ = s 
{\text{\rm Per}}_s(E;\Omega) -s
{\text{\rm Per}}_s(E_1;\Omega)  
-  s \int_{\Omega\setminus E_1} \int_{E_2 \cap B_R} 
\frac{1}{|x-y|^{n+s}}\, dx \, dy   \\
&  \qquad \qquad \ +  s \int_{ E_1} \int_{E_2 \cap B_R} \frac{1}{|x-y|^{n+s}}\, dx \, dy.  
\end{split}\end{equation}
By using \eqref{700} with $F:=\Omega\setminus E_1$ and $F:= E_1$
(which have finite $s_0$-perimeter in $\Omega$, recall the footnote
on page \pageref{foot}), 
we have that the last two terms in \eqref{10} converge to zero as 
$s\searrow 0$, thus
\begin{eqnarray}\label{800}
\lim_{s\searrow0} I(s,R)=\lim_{s\searrow 0}\Big(s 
{\text{\rm Per}}_s(E;\Omega) -s
{\text{\rm Per}}_s(E_1;\Omega)\Big).
\end{eqnarray}
We now recall the notation in~\eqref{AS} and we write 
\begin{eqnarray*}
\alpha_s(E)\,|\Omega\setminus E_1| 
&=& s \int_{\Omega\setminus E_1} 
\int_{E_2 \cap (\CC B_R)} \frac{1}{|x-y|^{n+s}}\, dx \, dy   
\\
&&\ +\, \alpha_s(E)\,|\Omega\setminus E_1|
-s \int_{\Omega\setminus E_1} \int_{E_2 \cap (\CC B_R)} 
\frac{1}{|x-y|^{n+s}} \, dx \, dy,  
\end{eqnarray*}
and 
\begin{eqnarray*}
\alpha_s(E)\,|E_1|
&=& s \int_{ E_1}
\int_{E_2 \cap (\CC B_R)} \frac{1}{|x-y|^{n+s}}\, dx \, dy              
\\ &&+\,\alpha_s(E)\,|E_1|
-s \int_{E_1} \int_{E_2 \cap (\CC B_R)} 
\frac{1}{|x-y|^{n+s}} \, dx \, dy.
\end{eqnarray*}
By subtracting term by term, we obtain that
\begin{eqnarray*}
&& \alpha_s(E)\,\Big( |\Omega\setminus E_1|-|E_1|\Big)\\
&& \qquad\qquad= I(s,R)+
\left(
\alpha_s(E)\,|\Omega\setminus E_1|
-s \int_{\Omega\setminus E_1} \int_{E_2 \cap (\CC B_R)}
\frac{1}{|x-y|^{n+s}} \, dx \, dy\right)
\\ &&\qquad\qquad\quad -\,\left(
\alpha_s(E)\,|E_1|
-s \int_{E_1} \int_{E_2 \cap (\CC B_R)}
\frac{1}{|x-y|^{n+s}} \, dx \, dy\right)
.\end{eqnarray*}
As a consequence, by using~\eqref{XX00}
(applied here both with~$F:=\Omega\setminus E_1$ and~$F:=E_1$),
\begin{equation}
\label{801}
\lim_{R\rightarrow+\infty}\lim_{s\searrow0}
\left[ \alpha_s(E)\,\Big( |\Omega\setminus E_1|-|E_1|\Big)
-I(s,R)\right]=0.
\end{equation}
Now, if  $|\Omega\setminus E|=|E\cap \Omega|$ then $|\Omega\setminus E_1|-|E_1|=0$,
and from \eqref{800}, \eqref{801}, and Corollary \ref{Cor} we get
$$
0=\lim_{R\rightarrow+\infty}\lim_{s\searrow0} I(s,R)=\lim_{s\searrow 0}s 
{\text{\rm Per}}_s(E;\Omega) -\mathcal M(E\cap \Omega),
$$
which proves that $E \in \mathcal E$ and $\mu(E)=\mathcal M(E\cap \Omega)$.
This establishes Theorem~\ref{TF1}(i).

On the other hand, if $|\Omega\setminus E|\neq|E\cap \Omega|$, then
by~\eqref{800}, \eqref{801}, and Corollary \ref{Cor}
we obtain the existence
of the limit
\begin{eqnarray*}
&&\Big( |\Omega\setminus E_1|-|E_1|\Big)\, \lim_{s\searrow0}
\alpha_s(E) \\
&& \qquad \qquad \qquad =
\lim_{R\rightarrow+\infty}\lim_{s\searrow0}
\alpha_s(E)\,\Big( |\Omega\setminus E_1|-|E_1|\Big)
\\[1ex]
&& \qquad \qquad  \qquad =
\lim_{R\rightarrow+\infty}\lim_{s\searrow0}\bigg\{ 
\left[ \alpha_s(E)\,\Big( |\Omega\setminus E_1|-|E_1|\Big)
-I(s,R)\right]
+\, I(s,R)\bigg\}
\\[2ex]
&& \qquad  \qquad \qquad = \mu(E)-\mu(E_1)= \mu(E)-\mathcal{M}(E\cap\Omega),\end{eqnarray*}
which completes the proof of
Theorem~\ref{TF1}(ii).~\hfill$\square$

\subsection{Construction of Example \ref{EXX}}

We start with some preliminary computations.
Let~$a_k:=10^{k^2}$, for any~$k\in\N$,
and let
$$
I_j:=\bigcup_{k\in\N} \big[a_{4k+j}, a_{4k+j+1}\big), \quad \text{for} \ j=0,1,2,3.
$$
Notice that~$[1,+\infty)$ may be written as the disjoint union
of the~$I_j$'s.
Let~$\varphi\in C^\infty\big([0,+\infty),\,[0,1]\big)$ be such that
$\varphi=0$ in~$[0,1]\cup I_0$,
$\varphi=1$ in~$I_2$, and then
$\varphi$  smoothly interpolates between~$0$ and~$1$
in~$I_1\cup I_3$.

We claim that 
there
exist two sequences~$\nu_{0,k}\rightarrow+\infty$ and~$\nu_{1,k}\rightarrow+\infty$
such that
\begin{equation}\label{DUE}
\lim_{k\rightarrow+\infty}\int_0^{+\infty} \varphi(\nu_{0,k} x) e^{-x}\,dx
=0\; \ {\mbox{ and }}\;
\lim_{k\rightarrow+\infty}\int_0^{+\infty} \varphi(\nu_{1,k} x) e^{-x}\,dx
=1.
\end{equation}
To check~\eqref{DUE}, we take~$\nu_{0,k}:=a_{4k+1}/k$
and~$\nu_{1,k}:=a_{4k+3}/k$. We observe that, by 
construction,~$\varphi=0$
in~$\big[a_{4k}, a_{4k+1}\big)$ and~$\varphi=1$ in~$\big[a_{4k+2}, a_{4k+3}\big)$,
so~$
\varphi(\nu_{0,k}x)=0$
for any~$x\in [kb_{0,k}, \,k)$
and~$\varphi(\nu_{1,k}x)=1$ in~$[kb_{1,k},\,k)$, where
$$ b_{0,k}:=\frac{a_{4k}}{a_{4k+1}}=10^{-(8k+1)}
\;{\mbox{ and }}
b_{1,k}:=\frac{a_{4k+2}}{a_{4k+3}}=10^{-(8k+5)}.$$
We deduce that
\begin{eqnarray*}
&& \int_0^{+\infty} \varphi(\nu_{0,k} x) e^{-x}\,dx\le
\int_0^{kb_{0,k}} e^{-x}\,dx+
\int_{k}^{+\infty} e^{-x}\,dx=1-e^{-kb_{0,k}}+e^{-k}
\\ {\mbox{and }}&&
\int_0^{+\infty} \varphi(\nu_{1,k} x) e^{-x}\,dx\ge
\int_{kb_{1,k}}^{k} e^{-x}\,dx=e^{-kb_{1,k}}-e^{-k}.
\end{eqnarray*}
This implies~\eqref{DUE} by noticing that
$$\lim_{k\rightarrow+\infty} k b_{0,k}=0=
\lim_{k\rightarrow+\infty} k b_{1,k}.$$

Now we construct our example by using the above function~$\varphi$
and~\eqref{DUE}.
We take~$\Omega := B_{1/2}$ and 
$E := \big\{ x =\left( \rho\cos\gamma , \rho\sin\gamma\right),  
\rho >1, \gamma\in\left[0,\theta\left(\rho\right) \right] \big\}\subset \R^2$, 
where $\theta\left( \rho\right) := \varphi \left(\log\rho\right)$. 

First of all, since $\Omega=B_{1/2}$ and $E \subset \R^n\setminus B_1$, it is easy to see that
$$
{\text{\rm Per}}_s(E;\Omega)  = \int_\Omega \int_E\frac{1}{|x-y|^{n+s}}\,dx\,dy 
\leq |\Omega| \, \int_{\R^n\setminus B_1} \frac{2^{n+s}}{|z|^{n+s}}\,dz<\infty
$$
for any $s \in (0,1)$ 
(notice that, since $E$ has smooth boundary, the fact that $E$ has
finite $s$-perimeter is also a consequence of Lemma 11 in \cite{CV11}). 
Then, recalling~\eqref{AS} we have
$$ \alpha_s(E) = 
 s \int_{1}^{+\infty} \int_{0}^{\theta\left(\rho\right)} \frac{\rho^{n-1}}{\rho^{n+s}} \, d\theta \, d\rho 
 =   s \int_{1}^{+\infty} \theta\left(\rho\right) \frac{1}{\rho^{1+s}} \, d\rho .
$$
Therefore, by the change of variable $\log\rho = r$, we have
\begin{equation*}
\alpha_s(E) = s \int_{0}^{+\infty} \varphi\left( r\right) e^{-rs}\,dr,
\end{equation*}
and, by the further change $rs=x$, we have
\begin{equation*}
\alpha_s(E) = \int_{0}^{+\infty} \varphi\left( \frac{x}{s}\right) e^{-x}\,dx. 
\end{equation*}
If we set $\nu ={1}/{s}$, the limit in \eqref{Rj} becomes the 
following:
\begin{equation*}
 \alpha(E) = \lim_{\nu\rightarrow \infty} \int_{0}^{+\infty} \varphi\left( \nu x\right) e^{-x}\,dx,
\end{equation*} 
and, by \eqref{DUE}, we get that such a limit
does not exist. 
This shows that the limit in \eqref{Rj} does not exist. 
Since~$|\Omega\setminus E|=|B_{1/2}|>0=|E\cap\Omega|$, 
by Theorem~\ref{TF1}(ii), the limit in \eqref{RR} does not exist 
either.~\hfill$\square$

\vspace{1mm}

\subsection{Construction of Example \ref{EXX special}}

It is sufficient to modify Example \ref{EXX} inside $\Omega=B_{1/2}$
in such a way that~$|\Omega\setminus E|=|E\cap\Omega|$. 
Notice that, since the set $E$ has smooth boundary, then
it has finite $s$-perimeter for any $s \in (0,1)$ (see Lemma 11 in \cite{CV11}). 
Then \eqref{Rj} is not affected by this modification and so
the limit in \eqref{Rj} does not exist in this case too.
On the other hand, the limit in \eqref{RR} exists, thanks
to Theorem~\ref{TF1}(i).~\hfill$\square$

\vspace{1mm}

\subsection{Construction of Example \ref{EX2}}

We take a decreasing
sequence~$\beta_k$ such that~$\beta_k>0$
for any~$k\ge1$,
$$ M:=\sum_{k=1}^{+\infty} \beta_k<+\infty$$
but
\begin{equation}\label{infinity}
\sum_{k=1}^{+\infty} \beta_{2k}^{1-s}=+\infty \qquad \forall \, s\in(0,1).
\end{equation}
For instance, one can take $\beta_1:=
\displaystyle\frac{1}{\log^2 2}$ and
$\displaystyle
\beta_k:=\frac{1}{k\log^2 k}$ for any~$k\ge2$.

Now, we define
\begin{eqnarray*}
&& \Omega:=(0,M)\subset\R,\\
&& \sigma_m:=\sum_{k=1}^m \beta_k,
\\ && I_m:=(\sigma_m,\sigma_{m+1}),\\
{\mbox{and }}&&
E:=\bigcup_{j=1}^{+\infty} I_{2j}.\end{eqnarray*}
Notice that~$E\subset\Omega$ and 
\begin{eqnarray}\label{above}
{\text{\rm Per}}_s(E;\Omega) & = & L(E,\CC E) \nonumber
\\
& \ge & \sum_{j=1}^{+\infty} 
L(I_{2j},I_{2j+1}) 
 \ =\  \sum_{j=1}^{+\infty}\int_{\sigma_{2j}}^{\sigma_{2j+1}} \int_{\sigma_{2j+1}}^{\sigma_{2j+2}}
\frac{1}{|x-y|^{1+s}}\,dx\,dy.
\end{eqnarray}
An integral computation shows that if~$a<b<c$ then
$$ \int_{a}^{b} 
\int_{b}^{c}
\frac{1}{|x-y|^{1+s}}\,dx\,dy=\frac{1}{s(1-s)}\Big[
(c-b)^{1-s}+(b-a)^{1-s}-(c-a)^{1-s}\Big].$$
By plugging this into~\eqref{above}, we obtain
\begin{equation}\label{above2}\begin{split}
& s(1-s){\text{\rm Per}}_s(E;\Omega) 
\\ &\qquad \ge \sum_{j=1}^{+\infty} \Big[
(\sigma_{2j+2}-\sigma_{2j+1})^{1-s}+(\sigma_{2j+1}-\sigma_{2j})^{1-s}-(
\sigma_{2j+2}-\sigma_{2j})^{1-s}\Big] \\
 & \qquad=  \sum_{j=1}^{+\infty}
\beta_{2j+2}^{1-s}+\beta_{2j+1}^{1-s}-
(\beta_{2j+2}+\beta_{2j+1})^{1-s}.
\end{split}\end{equation}
Now we observe that the map~$[0,1)\ni t\mapsto (1+t)^{1-s}$ is concave, therefore
$$ (1+t)^{1-s}\le1+(1-s)t\le 1+(1-s) t^{1-s}$$
for any~$t\in[0,1)$, that is
$$ 1+t^{1-s}-(1+t)^{1-s}\ge s t^{1-s}.$$
By taking~$t:=\beta_{2j+2}/\beta_{2j+1}$
and then multiplying by~$\beta_{2j+1}^{1-s}$, we obtain
$$ 
\beta_{2j+1}^{1-s}+\beta_{2j+2}^{1-s}-(\beta_{2j+1}+\beta_{2j+2})^{1-s}\ge
s \beta_{2j+2}^{1-s}.$$
By plugging this into~\eqref{above2} and using \eqref{infinity}, we conclude that
$$ 
{\text{\rm Per}}_s (E;\Omega)
\, \ge\, \frac{1}{1-s}
\sum_{j=1}^{+\infty}\beta_{2j+2}^{1-s}=+\infty \qquad \forall\, s\in (0,1),
$$
as desired.\hfill$\square$

\vspace{2mm}

\end{document}